\documentclass[reqno 11pt]{article}

\usepackage{amsmath,latexsym,amssymb,
            amscd, verbatim, graphics}

\def\sqr#1#2{{\vcenter{\hrule height.#2pt
        \hbox{\vrule width.#2pt height#1pt \kern#1pt
                \vrule width.#2pt}
        \hrule height.#2pt}}}
\def\square{\mathchoice\sqr64\sqr64\sqr{4}3\sqr{3}3}
\def\QED{\hfill$\square$\break}
\def\demo{\noindent{\bf Proof: }}

\newtheorem{Theorem}{\sc Theorem}[section]

\newtheorem{Corollary}[Theorem]{\sc Corollary}

\newtheorem{Remark}[Theorem]{\sc Remark}

\newtheorem{Definition}[Theorem]{\sc Definition}

\begin{document}

\baselineskip=13pt

\pagestyle{empty}

\ \vspace{1.7in}

\noindent {\LARGE\bf A Theorem
of Eakin and Sathaye and \\
Green's Hyperplane Restriction Theorem}

\vspace{.25in}

\noindent GIULIO \ CAVIGLIA, \  Department of Mathematics, U.C.
Berkeley, USA. \\{\it E-mail}: {\tt caviglia@math.berkeley.edu}

\vspace{2.4cm}

\begin{abstract}
A Theorem of Eakin and Sathaye relates the number of generators of
a certain power of an ideal with the existence of a distinguished
reduction for that ideal. We prove how this result can be obtained
as a special case of Green's General Hyperplane Restriction
Theorem.
\end{abstract}

\section{Introduction \hfill\break}
The purpose of these notes is to show how the following Theorem
\ref{EaSa}, due to Eakin and
Sathaye, can be viewed, after some standard reductions, as a corollary of
Green's General Hyperplane Restriction Theorem.\\

\noindent{\em {\sc Theorem \ref{EaSa}[Eakin-Sathaye]} Let $(R,m)$ be a
quasi-local ring with infinite residue field. Let $I$ be an ideal of $R$. Let $n$ and $r$ be
positive integers. If the number of minimal generators of
$I^i$, denoted by $v(I^i)$, satisfies
\[ v(I^i)<\binom{i+r}{r},\]
then there are elements $h_1,\dots,h_r$ in $I$ such that
$I^i=(h_1,\dots,h_r)I^{i-1}.$
}\\

Before proving Theorem \ref{EaSa} we have to recall some general facts
about Macaulay representation
of integer numbers. This is needed for the understanding of Green's Hyperplane Restriction Theorem. For more
details on those topics we refer the reader to \cite{G1} and \cite{G2}.

\subsection {Macaulay representation of integer numbers}
Let $d$ be a positive integer. Any positive integer $c$ can then be uniquely expressed as
\[
c=\binom{k_d}{d}+\binom{k_{d-1}}{d-1}+\dots + \binom{k_1}{1},
\]
where the $k_i$'s are non-negative and strictly increasing i.e
$k_d>k_{d-1}>\dots>k_1 \geq 0$. This way of writing $c$ is called the
$d$'th \emph{Macaulay representation} of $c,$ and the $k_i$'s are
called the $d$'th \emph{Macaulay coefficients} of $c.$
For instance, setting $c=13$ and $d=3$ we get $13=\binom{5}{3}+\binom{3}{2}+\binom{0}{1}.$

\begin{Remark} \label{Mac} {\rm
An important property of Macaulay representation is that the usual
order on the integers corresponds to the lexicographical order on
the arrays of Macaulay coefficients. In other words, given two positive
integer $c_1=\binom{k_d}{d}+\binom{k_{d-1}}{d-1}+\dots +
\binom{k_1}{1}$ and $c_2=\binom{h_d}{d}+\binom{h_{d-1}}{d-1}+\dots +
\binom{h_1}{1}$ we have $c_1<c_2$ if and only if
$(k_d,k_{d-1},\dots,k_1)$ is smaller lexicographically than
$(h_d,h_{d-1},\dots,h_1)$.}
\end{Remark}
\begin{Definition}{\rm Let $c$ and $d$ be positive integers.
We define $c_{<d>}$ to be
\[
c_{<d>}=\binom{k_d-1}{d}+\binom{k_{d-1}-1}{d-1}+\dots + \binom{k_1-1}{1}
\]
where $k_d,\dots,k_1$ are $d$'th Macaulay coefficients of $c$. We use
the convention that $\binom{a}{b}=0$ whenever $a<b.$
}\end{Definition}
\begin{Remark}\label{Mac1}{\rm
It is easy to check that if $c_1\leq c_2$ then ${c_1}_{<d>}\leq
{c_2}_{<d>}.$ This property, as we see in the following, allows us
to iteratively apply Green's Theorem and prove Corollary \ref{iter}.
}\end{Remark}

\subsection{Green's General Hyperplane Restriction Theorem}
Let $R$ be a standard graded algebra over an infinite field $K$. We can
write $R$ as $K[X_1,\dots,X_n]/I$ where $I$ is an homogeneous
ideal. Given a generic linear form $L$ we will denote by
$R_L=K[X_1,\dots,X_{n-1}]/I_L$ the restriction of $R$ to the hyperplane given
by $L$. Note that since $L$ is generic we can write it as
$L=l_1X_1+\dots+l_nX_n$ where $l_n\not = 0$, therefore $I_L$ is defined as
\[
I_L=(P(X_1,\dots,X_{n-1},(L/l_n)-X_n) \vert P\in I).
\]

\noindent We will denote by $R_d$ the $d$'th graded component of $R$.
Mark Green proved the following Theorem.

\begin{Theorem}[Green's General Hyperplane Restriction Theorem]\label{Green}
Let $R$ be a standard graded algebra over an infinite field $K$, and let $L$ be a generic linear form
of $R.$ Setting $S$ to be $R_L$, we have
\[
\dim_k S_d\leq (\dim_K R_d)_{<d>}.
\]
\end{Theorem}

The General Hyperplane Restriction Theorem first appeared in
\cite{G2}, where it was proved with no assumption on the
characteristic of the base field $K.$

A different, and more combinatorial, proof can be found in \cite{G1}
where the characteristic zero assumption is a working hypothesis. A
person interested in reading this last proof can observe that the
arguments in \cite{G1} also work in positive characteristic with
a few minor changes.

A direct corollary of Green's Theorem is the following
\begin{Corollary}\label{iter}
Let $R$ be a standard graded algebra over an infinite field $K$, and
let $L_1,\dots,L_r$ be generic linear forms of $R.$
Let $\binom{k_d}{d}+\binom{k_{d-1}}{d-1}+\dots + \binom{k_1}{1}$ be
the Macaulay representation of $dim R_d,$ and define
$S=R/(L_1,\dots,L_r)$. Then
\[
dim_K S_d\leq \binom{k_d-r}{d}+\binom{k_{d-1}-r}{d-1}+\dots + \binom{k_1-r}{1}
\]
\end{Corollary}
\demo
Note that $R_L$ is isomorphic to  $R/(L)$ and by Theorem \ref{Green}
one deduces  $dim_K (R/(L))_d\leq
\binom{k_d-1}{d}+\binom{k_{d-1}-1}{d-1}+\dots + \binom{k_1-1}{1}.$ On
the other hand by Remark \ref{Mac1} we can apply Green's Theorem
again and obtain the result by induction.
\QED

\section{The Eakin-Sathaye Theorem \hfill\break}
We now prove Theorem \ref{EaSa}. First of all note that  since
$v(I^i)$ is finite, without loss of generality we can assume that $I$
is also finitely generated: in fact if $J\subseteq I$ is a finitely
generated ideal such that $J^i=I^i$, the result for $J$ implies
the one for $I.$
Moreover, by the use of Nakayama's Lemma, we can replace $I$ by the
homogeneous maximal ideal of the fiber cone $S=\bigoplus_{i\geq 0} I^i/mI^i.$ Note that
$S$ is a standard graded algebra finitely generated over the infinite
field $R/m=K.$

Theorem \ref{EaSa} can be rephrased as:
\begin{Theorem}[E-S] \label{EaSa}
Let $R$ be a standard graded algebra finitely
generated over an infinite field $K.$ Let $i$ and $r$ be positive integers such that
\[
\dim_K(R_i)<\binom{i+r}{r}.
\]
Then there exist homogeneous linear forms $h_1,\dots,h_r$ such that
$(R/(h_1,\dots,h_r))_i$ is equal to zero.
\end{Theorem}
\demo
First of all note that
$\dim_K R_i \leq \binom{i+r}{r}-1=\binom{i+r}{i}-1=\binom{i+r-1}{i}+\binom{i+r-2}{i-1}+\dots+\binom{i+r-j}{i-j+1}+\dots
+\binom{r}{1}.$ This can be proved directly or by using
Remark \ref{Mac}. In fact one can
first order the array of Macaulay coefficients using the
lexicographic order and then note that the
previous array of $(i+r,0,\dots,0)$ is given by
$(i+r-1,i+r-2,\dots,r).$

Let $L_1,\dots,L_r$ be generic linear forms.
By Corollary \ref{iter} we have
\[
\dim_K (R/(L_1,\dots,L_r))_i\leq
\binom{i-1}{i}+\binom{i-2}{i-1}+\dots
+\binom{0}{1}
\]
The term on the right hand side is zero and therefore the theorem is
proved. \QED

\end{document}